\documentclass[12pt,leqno]{amsart}
\usepackage[dvips]{graphicx}
\usepackage{amsfonts}
\usepackage{amssymb}
\usepackage{amsmath}
\usepackage{amscd}
\usepackage{graphicx}
\usepackage[utf8]{inputenc}   
\usepackage[T1]{fontenc}      
\usepackage{geometry}         
\usepackage[francais]{babel}  
\usepackage{cjhebrew}
\def\DATE{\today}
\newtheorem{theorem}{Theorem}
\newtheorem{definition}[theorem]{Definition}
\newtheorem{corollary}[theorem]{Corollary}

\newtheorem{lemma}[theorem]{Lemma}

\newtheorem{proposition}[theorem]{Proposition}

\newcommand\C{\mathbb{C}}
\newcommand\R{\mathbb{R}}

\newcommand\ra{\rangle}

\newcommand\N{\mathbb{N}}

\newcommand\h{\mathfrak{h}}

\newcommand\K{\mathbb{K}}
\newcommand\Z{\mathbb{Z}}
\newcommand\Q{\mathbb{Q}}

\newcommand\pf{\noindent{\it Proof. }}

\newcommand{\kk}{^*\K}

\newcommand\lr{\left\{ \begin{array}{l}}
\def\ds{\displaystyle}

\title{Perturbations of  polynomials and applications}
\author{Elisabeth Remm,}
\date{}
\address{ E.R: Universit\'{e} de Haute Alsace, IRIMAS,4 rue des Fr\`{e}res Lumi\`{e}%
re, 68093 Mulhouse }
\email{elisabeth.remm@uha.fr}

\begin{document}

\maketitle

\begin{center}
A la m\'emoire de Robert Lutz
\end{center}

\begin{abstract} After reconsidering the theorem of continuity of the roots of a polynomial in terms of its coefficients in the deformation framework, we study the stability of the greater common divisor of two polynomials compared to perturbations on their roots. We apply this results to the study of deformations of a linear operator in finite dimension and in particular to the roots study of deformed matrices.
\end{abstract}


\section{Introduction}


Consider a given complex polynomial. In many cases the study of the polynomial which are close to it can by of great interest for the approximation of the roots, the control of little perturbations on the problem datas. As an example we can consider an Input/Output system. They naturaly  appear in fluids mechanics, signal theory,...  
An Input/Output linear system is generally modeled by an differential equation of the following form:
\begin{equation}
\label{trans}
\ds a_0 x(t)+a_1\frac{dx(t)}{dt}+\cdots+a_n\frac{d^nx(t)}{dt^n}=b_0 y(t)+b_1\frac{dy(t)}{dt}+\cdots+b_p\frac{d^py(t)}{dt^p}.
\end{equation}
If the Input function is  $x(t)=x_0 \sin(\omega y+ \varphi)$, the Output function is defined from the ratio
$$\ds H(i \omega)=\frac{b_0 +b_1(i \omega)+\cdots+b_m(i \omega)^m}{a_0 +a_1(i \omega)+\cdots+a_n(i \omega)^n}.$$
This function is called the transfert function. Let us consider the Laplace transform of  (\ref{trans}). If $X(p)$ and $Y(p)$ respectively refer to the Laplace transforms of $x(t)$ and $y(t)$, the transfert function corresponds to rational function
$$H(p)=\frac{Y(p)}{X(p)}=\frac{b_0 +b_1p+\cdots+b_mp^m}{a_0 +a_1p+\cdots+a_np^n}.$$
In order to reduce this rational expression, we have to compute the GCD of $Y(p)$ and $X(p)$.

Thus, formally, the problem of integration goes through a classical algebraic problem on rational fractions. However, the data are often linked to uncertain issues related to either measuring devices or numerical mathematical analysis techniques. We are therefore faced with the study of polynomial perturbations, how the roots of a polynomial evolve in a polynomial approximation or how the GCD evolves in these approximation  problems. 
A classic result in this area is the continuity of the roots in relation to the coefficients. Many results are known in this field. For example, let’s mention the study of pseudo-zeros of a polynomial or the study of GCD (\cite{Ga1,Ga2}). 

One of the fundamental results in complex polynomial study is the continuity of roots in relation to the coefficients. This theorem has been several times discussed and even recently \cite {K.H}. We here propose another approach based on the Lutz and Goze  \cite{L.G}  definition of perturbations with coefficients in a valuation ring. This notion of perturbations, sometimes called deformations, is largely used in algebra in the study a certain algebra laws like deformations of associative algebras, Lie algebras, Jordan algebras \cite{Ge}. The perturbations  used  here are based on and follow  this approach. It should also be noted that the study of polynomials on a valuation rings has been the subject of a number of publications. Let us mention for example \cite{Stef1}  which concerns the factorization problem. We can also note a presentation by the same author \cite{Stef2} concerning a synthesis of several works on polynomials.

In this work, we consider these results in an algebraic context adapted to the perturbations where they are described in an extension of the basic field associated with a valuation ring.  That allows a natural definition of infinitely small elements which is the  cornerstone of the perturbation theory.
\tableofcontents

\section{The non-archimedean field $\R^*$}

\subsection{Valuation rings}

Recall that an integral domain is a nonzero commutative ring in which the product of any two nonzero elements is nonzero. In these rings we have the following fundamental property 
$$ {\rm If} \ a \neq 0, \ {\rm then} \  ab = ac \Longrightarrow  b = c.$$

If $A$ is an integral domain, the field of fractions $F_A$  of $A$ is the smallest field containing $A$ that is there is an injective ring homomorphism $A \ra F_A$ such that any injective ring homomorphism from $A$ to a field can be factorized through $F_A$. To construct it, we consider the equivalence relation on the product  $A \times A^*$ where $A^*=A-0$:
$$(a,b) \mathcal{R} (c,d) \Leftrightarrow ad=bc$$
and $F_A$ is the quotient set associated with this equivalence relation. The equivalence class of $(a,b)$ is denoted by $\ds \frac{a}{b}$.

\begin{definition}
A valuation ring $A$ is an integral domain such that for any element $x \in F_A$ the field of fractions of $A$, one of the two elements $x$ or $x^{-1}$ is in $A$ where $x^{-1}$ is the inverse of $x$ in $F_A$.
\end{definition}
Let $\K$ be a field. A subring $A$ of $\K$ which satisfies
$$\forall x \in \K \setminus \{0\},  \ x \ {\text{\rm or}} \ x^{-1} {\text{\rm \ is  in \ }} A$$
is a valuation ring. It is called a valuation ring of $\K$. In this case $\K=F_A$.  

\begin{proposition}
Any valuation ring is a local ring that is has a unique maximal (left and right) ideal.
\end{proposition}
Let  $A$ be a valuation ring. The maximal ideal $\frak{m}$  of $A$ is the set of non inversible elements of $A$. In particular $A/\frak{m}$ is a field.

\noindent{Examples.}
\begin{enumerate}
  \item Any field $\K$  is a valuation ring of $\K$. But this example has very little interest.
  \item Let $\K$ be a field of characteristic $0$. Let us consider the ring $\K[[X]]$ of formal power series in the variable $X$ with coefficients in $\K$ that is any element of $\K[[X]]$ writes  $\ds \sum_{n \geq 0} a_nX^n$ with $a_n \in \K$. This ring is an integral domain. Let 
  $\K((X))$ be its field of fractions. Its elements are the Laurent series in $X$ that is the series 
$\ds \sum_{n\geq -M}a_n X^n$ where $M \in \Z$. The ring $\K[[X]]$  is a valuation ring of $\K((X))$. Its maximal ideal is the set of formal series with no term in degree $0$, that is $\ds \sum_{n\geq 1}a_n X^n$.
\item  Let $\K$ be an ordered field. This means that there is a total ordering of its elements that is compatible with the field operations.  The characteristic of such field is $0$ and $\Q$ is its prime field. Thus we can consider that $\Z$ is a subring of $\K$. An element $x \in \K$ is called finite if there are two integers $m,n \in \Z$ such that $n < x <m$. An element which is not finite is called infinite. The subset $A$ of finite elements of $\K$ is a valuation ring of $\K$. Then it contains a maximal ideal $\mathcal{I}$ whose elements are called infinitesimals. Then any infinitesimal is finite. In conclusion, in an ordered field we have: \begin{enumerate}
  \item the elements of $A$: the finite elements,
  \item the elements of $\mathcal{I} \subset A$: the infinitesimals,
  \item the elements $x\notin A$ such that $x^{-1} \in \mathcal{I}$ (see the following proposition) : the infinitely large elements. 
\end{enumerate}
Any perturbation or deformation problem will be modeled using this type of fields. Moreover the basic rule in analysis which is that the inverse of an infinitely large is infinitesimal is given by the following property:
\end{enumerate}

\begin{proposition}
Let $A$ be a valuation ring of a field  $\K$ and let  $\mathfrak{m}$ be its maximal ideal. We have 
$$x \in \K\setminus A \Rightarrow  x^{-1} \in \mathfrak{m}.$$
\end{proposition}
However this model of infinitesimal numbers is incomplete. Consider for example the real analysis. To consider the field of reals numbers as an valuation ring on itself has no interest. Since $\R$ is an archimedean field, the only infinitesimal is zero. We must therefore consider real numbers as belonging to a larger valuation ring. In this case we will be able to speak of infinitesimal but such numbers will be not in $\R$. It is also necessary to make a link between the field of this valuation ring and $\R$. Such extension exists and is given by the valuation ring of Robinson. However  there exists another valuation ring that is used in the deformation theory of algebraic structures: the formal series. 
\subsection{The valuation ring of formal series}
We have introduced this ring in Example 2. Consider $A=\K[[X]]$. Its field of fractions is $\K((X))$ and its maximal ideal is 
$$\frak{m}=\left\{\ds \sum_{n \geq 1} a_nX^n, a_n \in \K\right\}.$$
Recall that a formal series $\sum_{n \geq 0} a_nX^n$ is a unit in $A$ if and only if $a_0 \neq 0.$ The map $v : A \rightarrow \K$ given by $v(\sum_{n \geq k} a_nX^n)=k$ with $a_k\neq 0$ is a valuation and the ideal generated by a formal series $f$ is constituted of the formal series $g$ such that $v(g) \geq v(f)$ this is equivalent to say that $g=X^{v(f)}h, h\in A$.

In this case any formal series in $\frak{m}$ will be called an infinitesimal and the elements of $\K[[X]]$ will be called finite. In the following section we will see an important application of this valuation ring in the theory of deformation of algebraic structures.

\subsection{The valuation ring of Robinson}

There is a valuation ring, built from the field of real numbers $ \R$ to find intuitive notions about infinitely small or large elements and which are formally defined as above. Let $\R$  be the ordered field of real numbers. Notice  that at this level the archimedean property  of this field implies  that all the elements are finite and that there is only one infinitely small element, the element $0$. To build infinitely small linked to $\R$, we consider an ordered field  extension of $\R$, denoted $ ^*\R$, which is a non archimedean field, called field of hypereals  (we will not build it here, but we will  characterize it further by its most fundamental properties). As this field is ordered, we can consider the  subring $A$ of finite elements. We have seen that it is a valuation ring of $ ^*\R$.
Its maximum $\mathcal{I}$ ideal is the ideal of infinitely small elements. As $\R$ is a subfield of $^*\R$, all elements of $\R$ are finite elements and for any positive infinitely small element  $\varepsilon \in  \, ^*\R$  there exists $a \in \R$  with 
$$\varepsilon < a.$$
In addition, as the residual field $A/ \mathcal{I}$ is isomorphic to $\R$, for any finite element $x \in \, ^*\R$, there is an unique real $a \in \R$ such that $x-a$ is in $\mathcal{I}$. This real $a$ depends only on $x$. It will be denoted $^\circ x$

One of the reasons to consider this valuation ring rather than another to study properties built on $\R$ (or $\C$) is because of the existence of a transfer principle that can be summarized as follows: Suppose a proposition is true in $ ^*\R$, a proposal that is expressed by functions of a finite number of variables and quantifiers. So such a proposal is still true in $ \R$, proposal obtained by considering only the universal quantifier on $\R$. Maybe, another reason to use this ring rather than the valuation ring of formal series is that it is easier and more natural to describe an infinitely small neighborhood of a given point.

\noindent{\bf Vocabulary.} Let us resume the vocabulary introduced in the general study of evaluation rings. We consider the field of real numbers $\R$ and the valued extension of Robinson $^*\R$. Let $\mathcal{I}$ the maximal ideal of the associated valuation ring.
\begin{enumerate}
  \item An element of $\mathcal{I}$ is called an infinitesimal.
  \item An element $x \in \, ^*\R$ such that $x^{-1} \in \mathcal{I}$ is called infinitely large.
  \item An element $x \in \, ^*\R$ such that $x^{-1} \in \mathcal{I}$ which is not  infinitely large is called finite.
\end{enumerate}

 Since $\R$ can be considered as a subfield of $^*\R$, the elements of $\R$ are finite. We deduce that for any $x \in \, ^*\R$ which is not infinitely large, there exists a unique element $^\circ x \in \R$ such that
 $$x=\, ^\circ x + \varepsilon$$
 with $\varepsilon \in \R$.  If $\K=\C$, we have a similar vocabulary and an element $z \in \, ^*\C$ is finite if and only if its module $|z|$ is a finite number in $^*\R$. 
 

\subsection{Deformation in multilinear algebra}

The notion of deformations in multilinear algebra has experienced important developments in recent times, the most notable being certainly that of associative algebra deformations and the theory of quantization deformations in mathematical physics. 

Recall that a $\K$-algebra is a $\K$-vector space equipped with  a multiplication $\mu$ that is a bilinear map on $V$. Let us assume that $V$ is finite dimensional and let  $\{e_1,\cdots,e_n\}$ be a fixed basis of $V$. A bilinear map on $V$ is given by its structure constants
$$\mu(e_i,e_j)=\sum_{k=1}^n C_{i,j}^k e_k$$
and we can identify $\mu$ with the point $(C_{i,j}^k)$ of $\K^{n^3}$. 
In many cases (associative algebras, anti-associative algebras, Lie algebras), $\mu$ satisfies a quadratic relation, for example in the associative case, we have the relation $\mu(\mu(X,Y),Z)=\mu(X,\mu(Y,Z))$ that is translated by an algebraic equation in $\K^{n^3}$. If we consider the valuation ring $\K[[t]]$ of $\K$, a formal deformation of $\mu$ is a point $(C_{i,j}^k (t))$ of $\K[[t]]^{n^3}$ such that $C_{i,j}^k (0)=C_{i,j}^k$ which satisfies the same algebraic relation as the point $C_{i,j}^k$. If $\mu_t$ corresponds to the formal multiplication associated with $(C_{i,j}^k (t))$, then we can write
$$\mu_t=\mu + \sum_{s \geq 1}\varphi_s$$
where $\varphi_{s}$ is the bilinear map associated to $(C_{i,j}^k (s))$, and $\sum_{s \geq 1}\varphi_s$ is a formal bilinear map with coefficients in the maximal ideal of $\K[[t]]$, that is with infinitesimal coefficients. A lot of works concerns the properties of the bilinear maps $\varphi_1$ and $\varphi_2$. We can refer for example  \cite{Ge,Ma-Def,G.A.Rig, G.R.Val}.

A similar and probably equivalent notion of deformation of a multilication of an algebra can be do using the valued ring of Robinson. Let $\mu$ be a bilinear map on the $\K$-vector space $V$ ($\K=\R$ or $\C$).  We consider now the field $^*\K$ and the $^*\K$-vector space $^*V$. We know (see \cite{L.G}) that the dimension of $^*V$ is equal to the dimension of $V$. Moreover, if $\{e_1,\cdots,e_n\}$ is a basis of $V$, since we can consider that any vector of $V$ is a finite vector of $^*V$, then $\{e_1,\cdots,e_n\}$ is also a basis of $^*V$. We deduce that $\mu$ considered as a bilinear map on $^*V$ writes also $\mu(e_i,e_j)=\sum_{k=1}^n C_{i,j}^k e_k$ with $C_{i,j}^k \in \K \subset  \, \!  ^*\K$. In this context we can define a notion of deformation in $^*\K$.

\begin{definition}
A $^*\K$-deformation (called also a perturbation) of the bilinear map $\mu$ on $V$ is a bilinear map $^*\mu$ on $^*V$ which satisfies
$$^*\mu(e_i,e_j)=\sum_{k=1}^n D_{i,j}^ke_k$$ with
$$^{\circ}(D_{i,j}^k)=C_{i,j}^k.$$
\end{definition}
Recall a result which can be found in \cite{Go-Per}.
\begin{proposition}
Let $(\alpha_1,\cdots,\alpha_p)$ be a vector of $(^*\K)^p$. There exist infinitesimals $\varepsilon_1,\cdots,\varepsilon_l$ and linearly independent vectors $v_1,\cdots,v_l$ with $l \leq p$ of $\K^n$ such that
$$(\alpha_1,\cdots,\alpha_p)=\varepsilon_1v_1+\varepsilon_1\varepsilon_2v_2+ \cdots+ \varepsilon_1\varepsilon_2\cdots\varepsilon_lv_l.$$
\end{proposition}
Applying this decomposition to a $^*\K$-deformation $^*\mu$ of a bilinear map $\mu$ on $V$ we obtain
$$ ^*\mu = \mu + \varepsilon_1\varphi_1+\varepsilon_1\varepsilon_2\varphi_2+ \cdots+ \varepsilon_1\varepsilon_2\cdots\varepsilon_l\varphi_l$$
where $\varphi_1,\cdots,\varphi_l$ are linearly independent bilinear maps on $V$.

Assume now that $\mu$ satisfies a quadratic identity $\mu \circ \mu =0$. For example, if $\mu$ is associative, this identity is 
$$\mu \circ \mu (X,Y,Z)=\mu(X,\mu(Y,Z))-\mu(\mu(X,Y),Z).$$ For any bilinear maps $\varphi_1$ and $\varphi_2$ we denote by
$\varphi_1 \bullet \varphi_2$  the trilinear map obtained by linearization of  $\circ$:  $\mu \bullet \mu= 2 \mu \circ \mu$. If $\mu_t$ is a formal deformation of $\mu$, then $\mu_t =\mu + \sum_{k\geq 1} t^k\varphi_k$ and $\mu_t \circ \mu_t=0$ is equivalent to the infinite system:
$$\sum _{i+j=2k+1}\varphi_i \bullet \varphi_j=0, \ \ \sum _{i+j=2k}\varphi_i \bullet \varphi_j + \varphi_k \circ \varphi_k=0$$
for any $k \geq 0$ with the notation $\varphi_0=\mu.$
We know the interpretation of this infinite linear system only for the degree one, there exists a cohomological complex which permits to say that the first equation $\varphi_0 \bullet \varphi_1=0$ is equivalent to say that $\varphi_1$ is a $2$-cocycle associated with this cohomology, and also for the degree $2$ and the second equation
$$\varphi_0 \bullet \varphi_2+\varphi_1 \circ \varphi_1=0$$
is the basis of the equations of quantization deformations. We have not interpretation for the degrees greater than $3$ except  that this infinite system has a solution as soon as the previous cohomology is trivial in degree $3$.

If $^*\mu$ is a $^*\K$-deformation, then 
$$ ^*\mu = \mu + \varepsilon_1\varphi_1+\varepsilon_1\varepsilon_2\varphi_2+ \cdots+ \varepsilon_1\varepsilon_2\cdots\varepsilon_l\varphi_l$$ 
can be reduced to
$$\mu \bullet \varphi_1=0$$
and 
$$\varepsilon_1\varphi_1\circ \varphi_1 + \varepsilon_2 \mu \bullet \varphi_2+\varepsilon_1\varepsilon_2\varphi_1 \bullet \varphi_2 + \cdots +\varepsilon_1\varepsilon_2^2\cdots\varepsilon_l^2\varphi_l \circ \varphi_l=0.$$
This equation is a linear equation in the space of trilinear maps on $V$ with coefficients in $^*\K$. To solve it we have to compare these coefficients. \begin{enumerate}
  \item If $\varepsilon_2/\varepsilon_1 \in \mathcal{I}$, then $\varphi_1\circ \varphi_1=0$. In this case  $\mu+\varepsilon_1\varphi_1$ is already a deformation of $\mu$. 
  \item If $\varepsilon_1/\varepsilon_2 \in \mathcal{I}$ then $\mu \bullet \varphi_2=0$ implying that $\varphi_1\circ \varphi_1=0$ and in this case also $\mu+\varepsilon_1\varphi_1$ is already a deformation of $\mu$. 
  \item If $\varepsilon_2=a\varepsilon_1+\varepsilon_1\varepsilon'_2$ with $a \in \K$, $a \neq 0$, and $\varepsilon'_2\in \mathcal{I}$ then
  $$ \varphi_1\circ \varphi_1 + a\mu \bullet \varphi_2=0.$$
  The system becomes 
  $$
  \left\{
  \begin{array}{l}
  \mu \bullet \varphi_1=0   \\
   \varphi_1\circ \varphi_1 + a\mu \bullet \varphi_2=0  
\end{array}
\right.
$$
so we find again the system studied  in quantization deformations.  
\end{enumerate} 


\subsection{Deformation in linear algebra}

If $A$ is a square matrix of $gl(n,\K)$, then any $\kk$-deformation of $A$ is a matrix $^*A$ of $gl(n,\kk)$ such that all the  coefficients of $^*A-A$ are in $\mathcal{I}$. The properties of $^*A$ in relation with the properties of $A$ can be presented in terms of invariant linear subspaces. In this case, we have to study in a first time the deformations of polynomials. 

\section{Deformations of polynomials}
In this section, we assume that $\K=\C$ and $^*\C$ is the valued Robinson extension of $\C$. A finite (respectively infinitely) element of $^*\C$ is an element of $^*\C$ whose module is finite in $^*\R$ (respectively in $\mathcal{I}$). Let $\C[X]$ be the polynomial ring in $X$ over $\C$. The polynomials of $\C[X]$ have the form 
$$P(X)=a_0+a_1X+\cdots+a_nX^n$$
with $a_i \in \C$. An element of the ring $^*\C[X]$ have the form 
$$\widetilde{P}(X)=\alpha_0+\alpha_1X+\cdots+\alpha_mX^m$$
with $\alpha_i \in  \, ^*\C$. If all the coefficients $\alpha_i$ are finite, then there exists a unique polynomial $^\circ \widetilde{P}(X)=a_0+a_1X+\cdots+a_mX^m$ in $\C[X]$ such that $^\circ \alpha_i=a_i$. It is clear that the degree of $^\circ \widetilde{P}(X)$ is smaller than or equal to the degree of $P(X)$.

\begin{definition}
Let $P(X)$  be a polynomial in $\C[X]$ of degree $n$. A $^*\C$-deformation of $P(X)$ is a polynomial $\widetilde{P}(X)$ of $^*\C[X]$ of degree $n$ such that the coefficients of the polynomial  $\widetilde{P}(X)-P(X)$ belong to $\mathcal{I}$.
\end{definition}

\subsection{On the theorem of continuity of the roots}
Since $\C$ is algebraically close, any polynomial of $\C[X]$ and by transfert principle any element $\widetilde{P}(X)$ of $^*\C[X]$ admits a root. Since $\widetilde{P}(X)$ is a $^*\C$-deformation of $P(X)$, any root $\lambda$ of $\widetilde{P}(X)$  decomposed in $\lambda= r+\varepsilon$ with $r \in \C$, $P(r)=0$ and $\varepsilon \in \mathcal{I}$. The above remarks can be  summarized as follows:
\begin{proposition}
Let $Q(X)=\alpha_0+\alpha_1X+\cdots+\alpha_nX^n$ be a polynomial of  $^*\C[X]$ of degree $n \in \N$ such that all the coefficients $\alpha_k$ are finite. Assume that $\alpha_n \notin \mathcal{I}$. Then all the roots of $Q(X)$ are finite and if  $\xi$ is a root of $Q(X)$, then $u= \, ^\circ \xi$ is a root of  $P(X)=a_0+a_1X+\cdots+a_nX^n$ where $a_i=\, ^\circ \alpha_i$, $i=1,\cdots,n$.
\end{proposition}
Let us notice also that if $\xi_{i_1},\cdots,\xi_{i_k}$ are  roots of $Q(X)$ and if 
 $$^\circ\xi_{i_1}=\cdots=\, ^\circ\xi_{i_k}=u \ {\rm and } \ \forall j>0, \, ^\circ\xi_{i_{k+j}}\neq u$$
then $u$ is a root of $P(X)$ of multiplicity $k$. 

The theorem of  continuity of the roots of a polynomial  in relation to the coefficients can be formulated as follows:

 \begin{theorem}
Let $P(X) \in \C[X]$ be a polynomial of degree $n$ and  $Q(X)$ a $^*\C$-deformation of $P(X)$. Let  $u_i$ be a root of multiplicity $r_i$ of $P(X)$. There exist  linear maps
 $$L_i: \, ^*\C[X] \rightarrow \C$$
which depend only of $P(X)$ such that for any root   $\xi$ of $Q(X)$ with $^\circ \xi=u_i$ we have
 $$(\xi-u_i)^{r_i}=L_i( Q(X)-P(X))+\varepsilon ||(Q-P)(X)||$$
where $\varepsilon$ is infinitesimal in $^*\C$.
 \end{theorem}
 \noindent {\it Proof.} {\cite{L.G} Let $\xi$ be a root of $Q(X)$. Then $u=\, ^\circ \xi$ is a root of $P(X)$. Let $r$ be the multiplicity of  $u$. Then we have
 $$P(\xi)=P(u+\varepsilon)=\ds\frac{\varepsilon^r}{r!}(P^{(r)}(u)+\zeta)$$
 and $\zeta \in \,^*\C$ is infinitesimal. By hypothesis $P^{(r)}(u) \neq 0$ implying that 
 $$\ds\frac{\varepsilon^r}{r!}=\frac{P(\xi)}{P^{(r)}(u)+\zeta}$$
that is
 $$\varepsilon^r=\ds\frac{P(\xi)r!}{P^{(r)}(u)+\zeta}=-\ds\frac{(Q(\xi)-P(\xi))r!}{P^{(r)}(u)+\zeta}.$$
 The polynomial $H(X)=Q(X)-P(X)$ is in $\mathcal{I}[X]$. The linear map
  $$L(H)=-\ds\frac{H(u)r!}{P^{(r)}(u)}$$
is suitable.
 
 \noindent{\bf Remark.} The previous theorem, known as the theorem of continuous dependence of the roots of a polynomial on its coefficients is generally written in the following form \cite{K.H}: let $P(X)$ be a polynomial of degree $n$ in $\C[X]$. Consider the natural action of the symmetric group $\Sigma_n$ on $\C^n$: 
  $$\sigma (a_1,\cdots,a_n)=(a_{\sigma (1)},\cdots,a_{\sigma(n)})$$
with $\sigma \in \Sigma_n$. In particular, if $(a_1,\cdots,a_n)$ are the roots of $P(X)$ and if we consider the equivalent class of this $n$-uple modulo this action, we obtain a map: 
 $$\Phi: \C_n[X] \rightarrow \C^n/\Sigma_n.$$
 The metric $\ds\inf_{\sigma \in S_n}|a_i-a_{\sigma(i)}|$ defines a topology on $\C^n/\Sigma_n$. Since
 $\C_n[X] $ is a normed vector space, the  theorem of continuous dependence of the roots of a polynomial on its coefficients  said that  $\Phi$ is a continuous map. 

\subsection{The pseudozeros of a polynomial}


Let $P(X)$ be a polynomial of $\C[X]$ of degree $n$ and $\widetilde{P}(X) \in \, ^*\C[X]$ a $^*\C$-deformation of $P(X)$.

\begin{definition} We call pseudozero of $P(X)$ a root of a $^*\C$-deformation of $P(X)$.

\end{definition}

From the remarks of the previous section a pseudozero of $P(X)$ is a finite  element $\xi \in ^*\C$ such that $ ^\circ \xi=u$ is a root of $P(X)$. If $u_1,\cdots,u_k$ are the distinct roots of $P(X)$ and $\h (u_i)$  the set of the elements $\xi_i$ of $^*\C$ such that $\xi_i-u_i$ belongs to  $\mathcal{I}$,  we have

\begin{proposition}
Let $P(X)$ be a polynomial of $\C[X]$ of degree $n$ and $u_1,\cdots,u_k$ its roots. The set of pseudozeros of $P(X)$ is $\bigcup_{i=1}^k \h (u_i)$. 
\end{proposition}

 \medskip
 
 \noindent{\bf Example.} Let $P(X)=a_0+a_1X$ be a polynomial of degree $1,$  that is  $a_1 \neq 0$. Its root is $u_1=-a_0/a_1$. Let us consider a $^*\C$-deformation of $P(X)$:
 $$Q(X)=(a_0+\varepsilon_0)+(a_1+\varepsilon_1)X.$$ Its root is $\xi_1=-\frac{a_0+\varepsilon_0}{a_1+\varepsilon_1}$ that is
 $\xi_1=u_1+\varepsilon$. The goal is to evaluate $\varepsilon$ with respect to $\varepsilon_0$ and $\varepsilon_1$. We have
 $$\xi_1=u_1+\varepsilon=-\frac{a_0+\varepsilon_0}{a_1+\varepsilon_1}$$
 giving
 $$\varepsilon=\frac{a_0\varepsilon_1-a_1\varepsilon_0}{a_1(a_1+\varepsilon_1)}.$$
 \begin{lemma}(Goze decomposition \cite{Go-Per, L.G})
 Let $(\epsilon_1,\epsilon_2) \in \mathcal{I}$ be two infinitesimals. There exists two linearly independent vectors $(U_1,U_2)$ in $\C^2$ such that
 $$(\epsilon_1,\epsilon_2)=\alpha_1U_1+\alpha_1\alpha_2U_2$$
 with $\alpha_1,\alpha_2$ infinitesimals in $^*\C$.
 \end{lemma}
 If we moreover assume that the frame $(U_1,U_2)$ is orthonormal, then it is unique. Applying this decomposition to the pair $(\varepsilon_0,\varepsilon_1)$, we obtain
 $$(\varepsilon_0,\varepsilon_1)=\alpha_0U_0+\alpha_0\alpha_1U_1=\alpha_0(v_0,v_1)+\alpha_0\alpha_1(w_0,w_1)$$
 with $v_0w_1-v_1w_0 \neq 0$. Notice that it also includes the case corresponding to $\varepsilon_1=k\varepsilon_0$ which corresponds to $\alpha_1=0.$ This decomposition implies
 $$a_0\varepsilon_1-a_1\varepsilon_0=a_0(\alpha_0v_1+\alpha_0\alpha_1w_1)-a_1(\alpha_0v_0+\alpha_0\alpha_1w_0)=\alpha_0(a_0v_1-a_1v_0)+\alpha_0\alpha_1(a_0w_1-a_1w_0)$$
 
 First case: $a_0v_1-a_1v_0\neq 0$ that is the vector $(a_0,a_1)$ which corresponds to the polynomial $P(X)$ and the vector $U_0$ are linearly independent.  In this case we have
 $$\varepsilon \simeq \alpha_0 \frac{a_0v_1-a_1v_0}{a_1^2}.$$
 
 Second case: $a_0v_1-a_1v_0 = 0$ that is the vector $(a_0,a_1)$ which corresponds to the polynomial $P(X)$ and the vector $U_0$ are linearly dependent. Necessarily  $a_0w_1-a_1w_0 \neq 0.$ Then 
 $$a_0\varepsilon_1-a_1\varepsilon_0=\alpha_0\alpha_1(a_0w_1-a_1w_0)$$
 and 
 $$\varepsilon \simeq \alpha_0\alpha_1\frac{a_0w_1-a_1w_0}{a_1^2}.$$
 
\medskip

The technical tool used in this approximation is the Goze decomposition. In the general case, this decomposition has the following form:
\begin{lemma}
Let $(\varepsilon_1,\cdots,\varepsilon_n)$ a vector of $\mathcal{I}^n$. Then there exit infinitesimals $\alpha_1,\cdots,\alpha_n$ such that
$$
(\varepsilon_1,\cdots,\varepsilon_n)=\alpha_1 U_1+\alpha_1\alpha_2U_2+\cdots+\alpha_1\alpha_2\cdots \alpha_nU_n$$
with linearly independent vectors $U_1,\cdots,U_n$ in $\C^n.$
\end{lemma}
Let us notice that we have same decomposition in any valued extension of $\K$ (\cite{G.R.Val}) and in this paper we also discuss the unicity of such decomposition. 

For any vector $U=(u_0,\cdots,u_{n}) \in \C^{n+1}$, we denote $U(X)$ the polynomial $U(X)=u_0X+\cdots+u_{n}X^n$. Let $P(X)=a_0+a_1X+\cdots+a_nX^n$ be a polynomial of $\C^n[X]$ and $Q(X)$ a $^*\C$-deformation of $P(X)$ then
$$Q(X)=(a_0+\varepsilon_0)+(a_1+\varepsilon_1)X+\cdots+(a_n+\varepsilon_n)X^n$$ with $\varepsilon_i\in \mathcal{I}$ and 
$$(\varepsilon_0,\varepsilon_1,\cdots,\varepsilon_n)=\alpha_0U_0+\alpha_0\alpha_1U_1+\cdots + \alpha_0\alpha_1\cdots \alpha_nU_n$$
with linearly independent  vectors $(U_0,\cdots,U_n)$ in $\C^{n+1}$. We denote  $\Xi$ the vector $(\varepsilon_0,\cdots,\varepsilon_n)$ of $^*\C^{n+1}$. 

\noindent Notations: If $V=(v_0,v_1,\cdots,v_n)$ is a vector of $^*\K^{n+1}$, then $V(X)$ is the polynomial of $^*\K[X]$ given by
$V(X)=v_0+v_1X+\cdots+v_nX^n$. With these notations we can write
$$Q(X)=P(X)+\Xi(X)=P(X)+\alpha_0U_0(X)+\alpha_0\alpha_1U_1(X)+\cdots+\alpha_0\alpha_1\cdots\alpha_nU_n(X)$$
with linearly independent $(U_0(X),\cdots,U_n(X))$  in $\C^{n+1}[X].$ Let $\nu$ be a root of $Q(X)$. There is a root $u$ of $P(X)$ such that $\nu=u+\xi$ with $\xi \in \mathcal{I}$. Then
$$Q(\nu)=0=P(u+\xi)+\Xi(u+\xi).$$
But
$$P(u+\xi)=\xi P'(u)+\frac{\xi^2}{2!}P"(u)+\cdots +\frac{\xi^n}{n!}P^{(n)}(u)$$
and
$$\Xi(u+\xi)=\Xi(u)+\xi\Xi'(u)+\frac{\xi^2}{2!}\Xi"(u)+\cdots +\frac{\xi^n}{n!}\Xi^{(n)}(u) $$
and $P(u+\xi)+\Xi(u+\xi)=0$ implies
$$\xi P'(u)+\xi^2 A+ \alpha_0U_0(u)+\alpha_0B=0$$
where $A$ is finite in $^*\C$ and $B$ is infinitesimal.

\noindent 1. Assume that $P'(u) \neq 0$, that is $u$ is a simple root of $P(X)$.
\begin{lemma} If  $P'(u) \neq 0$, then $\Xi(u) \neq 0$. 
\end{lemma}
\pf In fact if $\Xi(u)=0$, we have
$$\xi P'(u)+\xi^2 A+\xi\Xi'(u)+\frac{\xi^2}{2}\Xi"(u)+\xi^2 A' =0$$
where $A'$ is infinitesimal.
This implies
$$P'(u)+\Xi'(u)=0$$
that is, since $\Xi(X) \in \mathcal{I}[X]$,
$$P'(u)=\Xi'(u)=0$$
then $\Xi(u) \neq 0.$ 

\medskip 

 We deduce
$$\xi P'(u) \simeq \Xi(u)$$ and
$$\xi \simeq -\frac{\Xi(u)}{P'(u)}.$$
But
$$\Xi(u)=\alpha_0U_0(u)+\cdots +\alpha_0\cdots \alpha_n U_n(u)$$
if $j$ is the smallest integer such that $U_j(u) \neq 0$, we have
$$\xi \simeq -\alpha_0\cdots \alpha_j\frac{U_j(u)}{P'(u)}.$$

\noindent 2.  Assume that $P'(u)=0$ and $P"(u) \neq 0$.
\begin{lemma} If  $P'(u)=0$ and $P"(u) \neq 0$, then $\Xi(u)+\xi\Xi'(u) \neq 0$. 
\end{lemma}
\pf It is analogous to  previous lemma's proof . 

\medskip 

We deduce that $Q(u+\xi)$ is equivalent to $\Xi(u)+\xi\Xi'(u)+\xi^2\ds \frac{P"(u)}{2}.$ If $\Xi(u)=0$ and $\Xi'(u) \neq 0$, we obtain
$$\xi \simeq -2\frac{\Xi'(u)}{P"(u)}.$$
If $\Xi(u) \neq 0$ and $\Xi'(u)=0$ then
$$\xi^2 \simeq -2\frac{\Xi(u)}{P"(u)}.$$
In the other cases, 
$$\frac{\Xi(u)+\xi\Xi'(u)}{\xi ^2}\simeq \ds -\frac{P"(u)}{2}.$$
We have to evaluate $\xi$. We have
$$\Xi(u)+\xi\Xi'(u)+\xi^2\ds \frac{P"(u)}{2}=\Xi(u)+\xi\left(\Xi'(u)+\xi\ds \frac{P"(u)}{2}\right).$$
If $\ds \frac{\Xi'(u)}{\xi} \simeq 0$, then 
$\xi^2 \simeq \ds -2\frac{\Xi(u)}{P"(u)}$ and we have already meet this case. If $\ds \frac{\xi} {\Xi'(u)}\simeq 0$, then 
$$\xi \simeq \ds -\frac{\Xi(u)}{\Xi'(u)}$$
with $$\ds \frac{\Xi(u)}{\Xi'(u)^2} \simeq 0.$$
If $\ds \frac{\xi}{\Xi'(u)} \simeq \frac{1}{a} $ with $a \neq 0, a \in \K$, then
$$\xi^2 \simeq -\frac{2\Xi(u)}{2a+P"(u)}$$
but in this case we must have
$$\Xi(u) \simeq \Xi'(u)^2.$$

If $\Xi(u)=\sum_{k=0}^n \alpha_0\alpha_1\cdots \alpha_k U_k(u)$,  $\Xi'(u)=\sum_{k=0}^n  \alpha_0\alpha_1\cdots \alpha_k U_k'(u)$ and if $j_0$ and $j_1$ denote the smallest integers such that $U_{j_0}(u) \neq 0$ and $U_{j_1}'(u) \neq 0$ then this last condition gives

$$\alpha_0\cdots \alpha_{j_0}U_{j_0}(u) \simeq  \alpha_0^2\cdots \alpha_{j_1}^2U'_{j_1}(u)^2$$

\noindent{\bf Remark.} We call rank of the vector $(\varepsilon_0,\cdots,\varepsilon_n)$ of $\mathcal{I}^{n+1}$ the greater index $k$ such that $(\varepsilon_0,\cdots,\varepsilon_n)=\sum_{i=1}^k\alpha_0\cdots\alpha_iU_i$ with $\alpha_i \neq 0$. If this rank is equal to $n+1$  then the polynomial $\Xi(X)$ have not root in $\K$. In particular $\Xi(u) \neq 0$. If $U_0(u)=0$, since the vectors $U_0,U_1,\cdots,U_n$ are linearly independent, then $u$ cannot be a root of all the polynomials $U_0(X),U_1(X),\cdots,U_n(X)$. 
\begin{theorem}\label{xi}
Let $P(X)$ be a polynomial of degree $n$ in $\C[X]$ and $Q(X)$ a $^*\C$-deformation of $P(X)$. Then $Q(X)=P(X)+\Xi(X)$ with
$$\Xi(X)=\varepsilon_0+\varepsilon_1 X+\cdots+\varepsilon_n X^n$$ with $\varepsilon_i\simeq 0$, that is $\varepsilon_i \in \mathcal{I}$. Considering the decomposition
$$(\varepsilon_0,\varepsilon_1,\cdots,\varepsilon_n)=\alpha_0U_0+\alpha_0\alpha_1U_1+\cdots + \alpha_0\alpha_1\cdots \alpha_nU_n$$
with $(U_0,\cdots,U_n)$ linearly independent in $\C^{n+1}$, if $u$ is a root of order $k$ of $P(X)$, then there exists $\xi \in \mathcal{I}$ such that $u+\xi$ is a root of $Q(X)$ and 
$$\xi^k \simeq  -\frac{\alpha_0\alpha_1\cdots \alpha_{j_0}U_{j_0}(u)}{P^{(k)}(u)}$$
as soon as $u$ is a root of the polynomilals $U_0(X),\cdots,U_{j_0-1}(X).$
\end{theorem}

 \noindent{\bf Example} Let us consider the matrix 
 $$M=\begin{pmatrix}
     1 &1    \\
     0 & 1  
\end{pmatrix}
$$
 Its characteristic polynomial is $P_M(X)=X^2-2X+1$ and $\lambda=1$ is a root of multiplicity $2$. Let us consider a $^*\C$-deformation (the deformations of matrices will be studied in more detail in Section \ref{Section 5}) given by the matrix
$$\widetilde{M}=\begin{pmatrix}
     1 &1    \\
   \varepsilon & 1  
\end{pmatrix}
$$ 
with $\varepsilon \in ^*\C$. Its characteristic polynomial is $P_{\widetilde{M}}(X)=X^2-2X+1-\varepsilon$. We have
$$P_{\widetilde{M}}(X)=P_M(X)+\Xi(X)$$
with $\Xi(X)=-\varepsilon.$ The corresponding vector of $^*\C^3$ is $(-\varepsilon,0,0)$ whose decomposition is trivial:
$$(\varepsilon,0,0)=\alpha_0U_0$$ with $\alpha_0=-\varepsilon$ and $U_0=(1,0,0)$ which corresponds to the polynomial $U_0(X)=-1$. Let $\lambda + \xi$ be a root of $P_{\widetilde{M}}$. Since $U_0(\lambda) \neq 0$, we have
$\ds\xi^2 \simeq -\frac{\alpha_0U_0(\lambda)}{P_M''(\lambda)}$ that is
$$\xi^2 \simeq \frac{\varepsilon}{2}.$$

\subsection{Linear deformation}
Let be $P(X) \in \C^n[X]$ and $Q(X)$ a $^*\C$-deformation of $Q(X)$ (or a formal deformation). We will say that $Q(X)$ is a linear deformation if $Q(X)=P(X)+\varepsilon R(X)$ where $R(X) \in \C^n[X]$. In other words $\Xi[X]=\varepsilon \left(b_0+b_1X+\cdots +b_nX^n\right)$ with $b_i \in \C$.
In this case, the associated decomposition of $\Xi[X]$ is $\Xi[X]=\alpha_0U_0(X)$ with $\alpha_0=\varepsilon$ and $U_0=(b_0,\cdots,b_n)$  the vector of $\C^{n+1}$ associated with $R(X)$. If $u+\xi$ is a root of $Q(X)$ with $P(u)=0$, then 
$$\xi \simeq -\varepsilon \frac{U_0(u)}{P'(u)}$$
if $u$ is a simple root of $P(X)$. If $u$ is a root of multiplicity $k$, then
$$\xi^k \simeq -\varepsilon \frac{U_0(u)}{P^{(k)}(u)}$$
 \section{GCD of perturbations}
 \subsection{Quasi-Zeros of a polynomial}
Let $P(X)$ be a polynomial of $\K[X]$ where $\K=\R$ or $\C$. We call quasi-zero of $P(X)$ a root of a perturbation $Q(X)$ of $^*\K[X]$. When the notion of perturbation is related with a notion of norm in the space of polynomials, the determination of quasi-zero is often difficult to describe. When the perturbation is defined from a valued ring, this determination is very easy. Consider $u$ a root of $P(X)$ and $h(u) \subset \; \! ^*\K$  the set $h(u)=\{u+\xi, \ \xi \in I\}$. The union of the $h(u)$ for all the roots $u$ of $P(X)$ is the set of the quasi-zeros of $P(X)$.

A natural problem then arises. Given two polynomials, how does their GCD behave when we perturb these two polynomials.
\subsection{The GCD of perturbations of polynomials}
Let $P_1(X)$ and $P_2(X)$ be two polynomials of $\K[X]$ and $D(X)$ their $GCD$. Let $Q_1(X)$ and $Q_2(X)$ be polynomials of $^*\K[X]$ which are perturbations respectively of $P_1(X)$ and $P_2(X)$. Usually, the GCD of $Q_1(X)$ and $Q_2(X)$ is not a perturbation of $D(X)$. 
\begin{lemma}
Let $Q_1(X)=Q_2(X)Q_3(X)+S(X)$ be the Euclidean division of $Q_1(X)$ by $Q_2(X)$ and $P_1(X)=P_2(X)P_3(X)+R(X)$ the Euclidean division of $P_1(X)$ by $P_2(X)$. Then $Q_3(X)$ and $S(X)$ are $^*\K$-perturbations of $P_3(X)$ and $R(X)$.
\end{lemma}
\pf If $P_1(X)=a_nX^n+ \cdots+a_0, \ P_2(X)=b_pX^p+ \cdots +b_0$, with $a_nb_p \neq 0$ and $n \geq p$, then $Q_1(X)=(a_n+\varepsilon_n)X^n+ \cdots+(a_0+\varepsilon_0), \ Q_2(X)=(b_p+\zeta_p)X^p+ \cdots +(b_0+\zeta_0)$ with $\varepsilon_i,\zeta_j \in \mathcal{I}$. Then the first coefficient of $Q_3(X)$ is $\frac{a_n+\varepsilon_n}{b_p+\zeta_p}$ which is of the form $\frac{a_n}{b_p}+\xi_{n-p}$ with $\xi_{n-p}\in \mathcal{I}$. Then 
$$\begin{array}{ll}
Q_1(X)& =Q_2(X)\left(\frac{a_n}{b_p}+\xi_{n}\right)X^{n-p} +\sum_{i=1}^p\left( \frac{a_{n-i}b_p-a_nb_{p-i}}{b_p} +\xi_{n-i}\right)X^{n-i}\\
& +\sum_{i=p+1}^n\left( a_{n-i}+\xi_{n-i}\right)X^{n-i}
\end{array}$$

with $\xi_{n-i} \in \mathcal{I}$. We see that the first remainder is a perturbation of the first remainder in the division of $P_1(X)$ by $P_2(X)$. By induction we obtain the proof.

 Using the Euclidean algorithm, the GCD of $P_1(X)$ and $P_2(X)$ is the last non nul remainder. It is also the case for $Q_1(X)$ and $Q_2(X)$. But this last remainder can be a singular perturbation (that is without condition on the degree of the perturbation) of the null polynomial, that is this last remainder can be a polynomial in $\mathcal{I}[X]$. 
 
 \begin{definition}
 We call the PGCD (Perturbed Greater Commun Divisor) of $Q_1(X)$ and $Q_2(X)$ the last reminder in the Euclidean algorithm  wich has not only infinitesimal coefficients  (that is the last remainder that is not in $\mathcal{I}[X]$).
 \end{definition}
From the previous calculus on the Euclidean division we can conclude: 
 \begin{theorem}
 Let $P_1(X)$ and $P_2(X)$ be two polynomials of $\K[X]$ and let $D(X)$ their $GCD$. Let $Q_1(X)$ and $Q_2(X)$ be polynomials of $^*\K[X]$ which are perturbations respectively of $P_1(X)$ and $P_2(X)$. Then the PGCD of $Q_1(X)$ and $Q_2(X)$ is a $^*\K$-perturbation of the GCD of $P_1(X)$ and $P_2(X)$.
 \end{theorem}
 We can summarize this result by
 $${\rm GCD}\left(P_1(X),P_2(X)\right)= \; ^\circ{\rm PGCD}\left(Q_1(X),Q_2(X)\right).$$
 
 \medskip
 
 \noindent{\bf Example.} If $P_1(X)=X^3-1$ and $P_2(X)=X^2-1$, then ${\rm GCD}(P_1(X),P_2(X))=X-1$. Let us consider the perturbations
 $Q_1(X)=X^3-\varepsilon_1X+(-1+\varepsilon_2)$ and $Q_2(X)=X^2+\varepsilon_3 X-1$. The Euclidean algorithm for $Q_1(X)$ and $Q_2(X)$ gives
 $$\begin{array}{l}
 X^3-\varepsilon_1X+(-1+\varepsilon_2)= (X^2+\varepsilon_3 X-1)(X-\varepsilon_3) +(1-\varepsilon_1+\varepsilon_3^2)X-(1-\varepsilon_2+\varepsilon_3)     \\  
 \\
 X^2+\varepsilon_3 X-1=((1-\varepsilon_1+\varepsilon_3^2)X-(1-\varepsilon_2+\varepsilon_3) )Q_3(X)+\ds\frac{\Xi(X)}{(1-\varepsilon_1+\varepsilon_3^2)^2}
\end{array}
$$
with $Q_3(X)=\ds\frac{1}{1-\varepsilon_1+\varepsilon_3^2}X+\frac{1-\varepsilon_2+2\varepsilon_3-\varepsilon_1\varepsilon_3+\varepsilon_3^3}{(1-\varepsilon_1+\varepsilon_3^2)^2}$ and 
$$\Xi(X)=2\varepsilon_1-2\varepsilon_2+3\varepsilon_3-\varepsilon_1^2+\varepsilon_2^2-\varepsilon_1\varepsilon_3-3\varepsilon_2\varepsilon_3+\varepsilon_3^3+\varepsilon_1\varepsilon_2\varepsilon_3+\varepsilon_1\varepsilon_3^2-\varepsilon_2\varepsilon_3^3.$$
Since $\Xi(X) \in \mathcal{I}(X)$, 
$${\rm PGCD}(Q_1(X),Q_2(X))=(1-\varepsilon_1+\varepsilon_3^2)X-(1-\varepsilon_2+\varepsilon_3).$$
This polynomial is a perturbation of ${\rm GCD}(P_1(X),P_2(X))=X-1.$
\subsection{Simplification of the transfert function}
In the introduction we recalled that in an I/O linear system the transfer function could be reduced to the form of a rational fraction
$$H(p)=\frac{Y(p)}{X(p)}=\frac{b_0 +b_1p+\cdots+b_mp^m}{a_0 +a_1p+\cdots+a_np^n}$$
If $D(p)={\rm GCD}(Y(p),X(p))$, then
$$H(p)=\frac{Y(p)}{X(p)}=\frac{D(p)Y_1(p)}{D(p)X_1(p)}=\frac{Y_1(p)}{X_1(p)}.$$
If the data of the linear system are known with some uncertainty, we have to reduce a perturbed transfert function
$$\widetilde{H}(p)=\frac{\widetilde{Y}(p)}{\widetilde{X}(p)}$$
and we can consider that $\widetilde{Y}(p),\widetilde{X}(p)$ are perturbations of $Y(p)$ and $X(p)$.  Let $\widetilde{D}(p)$ be the ${\rm PGCD}(\widetilde{Y}(p),\widetilde{X}(p))$. The Euclidean algorithm gives
$$
\begin{array}{l}
     \widetilde{Y}(p)=\widetilde{X}(p) Q_1(p)+R_1(p)  \\
      \widetilde{X}(p) =R_1(p)Q_2(p)+R_2(p)\\
      \cdots \\
      R_{k-1}(p)=R_k(p)Q_k(p)+\widetilde{D}(p)\\
     R_k(p)=\widetilde{D}(p) Q_{k+1}(p)+\Xi (p)
\end{array}
$$
with $\Xi (p) \in \mathcal{I}[X]$ and $R_j(X) \notin \mathcal{I}[X]$ for $j=1,\cdots,k.$ We deduce
$$\widetilde{X}(p)=D(p)\widetilde{X_1}(p)+\Xi_1(p), \ \ \widetilde{Y}(p)=D(p)\widetilde{Y_1}(p)+\Xi_2(p)$$
with $\Xi_1(p),\Xi_2(p) \in \mathcal{I}[X]$ and $\widetilde{X_1}(p),\widetilde{Y_1}(p)$ with finite coefficients in $^*\K$. Then
$$\widetilde{H}(p)=\frac{\widetilde{Y}(p)}{\widetilde{X}(p)}=\frac{D(p)\widetilde{Y_1}(p)+\Xi_2(p)}{D(p)\widetilde{X_1}(p)+\Xi_1(p)}=\frac{\widetilde{Y_1}(p)}{\widetilde{X_1}(p)}+\Xi(p).$$

\noindent{\bf Example.} Let us assume that, for example, 
$$H(p) = \frac{Y(p)}{X(p)}=\frac{p^3-1}{p^2-1}$$
and 
$$\widetilde{H}(p)=\frac{\widetilde{Y}(p)}{\widetilde{X}(p)}=\frac{p^3-\varepsilon_1p+(-1+\varepsilon_2)}{p^2+\varepsilon_3p-1}.$$
We have already computed the PGCD:
$${\rm PGCD}(\widetilde{Y}(X),\widetilde{X}(X))=\widetilde{D}(X)=(1-\varepsilon_1+\varepsilon_3^2)X-(1-\varepsilon_2+\varepsilon_3).$$
The Euclidean algorithm gives
$$
\begin{array}{l}
     \widetilde{Y}(p)=\widetilde{X}(p)(p-\varepsilon_3)+\widetilde{D}(p)    \\
    \widetilde{X}(p)=\widetilde{D}(p)Q_3(p)+\Xi_1 (p)  
\end{array}
$$
with $\Xi_1 (p) =\ds \frac{\Xi(p)}{(1-\varepsilon_1+\varepsilon_3^2)^2}\in \mathcal{I}[X].$

Then
$$\widetilde{Y}(p)=\widetilde{D}(p)\left((p-\varepsilon_3)Q_3(p)+1\right)+(p-3)\Xi_1 (p)$$
and
$$\widetilde{H}(p)=\frac{\widetilde{D}(p)((p-\varepsilon_3)Q_3(p)+1)+(p-3)\Xi_1 (p)}{\widetilde{D}(p)Q_3(p)+\Xi_1 (p)}.$$
Since $\widetilde{D}(p)=D(p)+\Xi_2(p)$ with $\Xi_2(X) \in \mathcal{I}[X]$, 
$$\widetilde{H}(p)=\frac{p^2+p+1}{p+1} +\Xi_3(p)$$
with $\Xi_3(X) \in \mathcal{I}[X].$ It may be interesting to evaluate the infinitesimal $\Xi_3(p)$. We have
$$\widetilde{H}(p)=p-\varepsilon_3+\frac{p(1-\varepsilon_1+\varepsilon_3^2)-1+\varepsilon_2-\varepsilon_3}{p^2+\varepsilon_3p-1}$$
and
$$\Xi_3(p)=\frac{p(1-\varepsilon_1+\varepsilon_3^2)-1+\varepsilon_2-\varepsilon_3}{p^2+\varepsilon_3p-1}-\frac{1}{p+1}-\varepsilon_3$$
After direct computations we obtain
$$\Xi_3(p)\simeq -\varepsilon_1\frac{p}{p^2-1}+\varepsilon_2\frac{1}{p^2-1}-\varepsilon_3\frac{p^3+p^2+p}{(p^2-1)(p+1)}$$
Considering the decomposition $(\varepsilon_1,\varepsilon_2,\varepsilon_3)=\alpha_1U_1+\alpha_1\alpha_2U_2+\alpha_1\alpha_2\alpha_3U_3$, we obtain 
$$\Xi_3(p)\simeq \alpha_1\left(-u_1\frac{p}{p^2-1}+v_1\frac{1}{p^2-1}-w_1\frac{p^3+p^2+p}{(p^2-1)(p+1)}\right)$$
where $U_1=(u_1,v_1,w_1)$ as soon as $-u_1p(p+1)+v_1(p+1)-w_1(p^3+p^2+p)\neq 0.$
The other cases can be treated similarly.
\section{Perturbations of matrices}{\label{Section 5}}
\subsection{Some generalities}
The vector space  $gl(n,\K)$ of square matrices of order $n$ is fibered by ne natural action of the Lie groups $GL(n,\K)$ of invertible matrices of order $n$, the orbit of a matrix $A$ is 
$$\mathcal{O}(A)=\{P^{-1}AP, \ P \in GL(n,\K)\}.$$
Some orbits are singular: $\mathcal{O}(A)=\{A\}$ as soon as $A=aId_n$  where $Id_n$ is the identity matrix,  In the other cases, the orbit is not reduced to a point. It is a provided with a differential manifold structure whose dimension is given by the dimension of its tangent space in $A$ which corresponds to
$$T_A(\mathcal{O}(A))=\{[M,A]=MA-AM, \ M \in gl(n,\K)\}$$

\noindent{\bf Example: $n=2$. $\K=\C$} 
\begin{enumerate}
  \item If $A= aId_2$, $\dim \mathcal{O}(A)=0.$
  \item If $A$ is diagonalizable with two distinct eigenvalues, then $\dim \mathcal{O}(A)=2.$
  \item If $A$ is not diagonalizable, then $\dim \mathcal{O}(A)=2.$
\end{enumerate}

\subsection{Characteristic polynomial of a perturbation of a square matrix}

Let $$A=\begin{pmatrix}
   a_{1,1}   & \cdots & a_{1,n}   \\
   a_{2,1}   &  \cdots & a_{2,n} \\
   \vdots & \cdots  & \vdots \\
   a_{n,1} & \cdots & a_{n,n}
\end{pmatrix}
$$
be a square matrix of order $n$. We denote by $q_{i_1,\cdots,i_k}^{(k)}(A)$ the determinant 
$$q_{i_1,\cdots,i_k}^{(k)}(A)=\det\begin{pmatrix}
  a_{i_1,i_1}    &a_{i_1,i_2}  &\cdots & a_{i_1,i_k}   \\
    a_{i_1,i_2}    &a_{i_2,i_2}  &\cdots & a_{i_2,i_k}   \\
  \vdots &  \vdots &  &  \vdots \\
  a_{i_k,i_1}  &a_{i_k,i_2}  &  \cdots & a_{i_k,i_k}
\end{pmatrix}
$$
with with $1\leq i_1<i_2<\cdots <i_k\leq n $ . The homogeneous polynomials
$$Q^{(k)}(A)=\sum_{1\leq i_1<i_2<\cdots <i_k\leq n}q_{i_1,\cdots,i_k}^{(k)}(A)$$
are the coefficients of the characteristic polynomial $C_A(X)$:
$$C_A(X)=X^n-Q^{(1}(A)X^{n-1}+Q^{(2)}(A)X^{n-2} + \cdots + (-1)^jQ^{(j)}(A)X^{n-j}+\cdots+(-1)^nQ^{(n)}(A).$$
Let us notice that
$$Q^{(1)}(A)={\rm tr}(A), \ \ {\rm and} \ \ Q^{(n)}(A)=\det (A).$$
Classically the polynomials $Q{(k)}(A)$ are written in terms of exterior products of $A$.  We associate to each polynomial $Q^{(k)}(A)$ the symmetric $k$-linear form $\Theta_{Q^{(k)}}$ on $gl(n,\K)$ of degree $k$ which satisfies
$$\Theta_{Q^{(k)}}(A,A,\cdots,A)=k!Q{(k)}(A).$$

\noindent{\bf Examples.}
\begin{enumerate}
  \item $n=2$. If 
  $$
  A=\begin{pmatrix}
    a_{1,1}  & a_{1,2}   \\
   a_{2,1}   &  a_{2,2} 
\end{pmatrix}
$$ then\begin{enumerate}
  \item $Q^{(1)}(A)=a_{1,1}+a_{2,2}, \ \Theta_{Q^{(1)}}(A)=Q^{(1)}(A).$
  \item $Q^{(2)}(A)=\det A=a_{1,1}a_{2,2}-a_{1,2}a_{2,1}, \ \Theta_{Q^{(2)}}(X,Y)=x_{1,1}y_{2,2}+x_{2,2}y_{1,1}-x_{1,2}y_{2,1}-x_{2,1}y_{1,2}$
\end{enumerate}
  \item $n=3$, $A=(a_{i,j})$
\begin{enumerate}
  \item $Q^{(1)}(A)={\rm tr}(A), \ \Theta_{Q^{(1)}}(A)=Q^{(1)}(A).$
  \item $Q^{(2)}(A)=a_{1,1}a_{2,2}-a_{2,1}a_{1,2}+a_{1,1}a_{3,3}-a_{3,1}a_{1,3}+a_{2,2}a_{3,3}-a_{3,2}a_{2,3}$
  
  $\Theta_{Q^{(2)}}(X,Y)= x_{1,1}y_{2,2}+y_{1,1}x_{2,2}-x_{2,1}y_{1,2}-y_{2,1}x_{1,2}+x_{1,1}y_{3,3}+y_{1,1}x_{3,3}-x_{3,1}y_{1,3}-y_{3,1}x_{1,3}+x_{2,2}y_{3,3}+y_{2,2}x_{3,3}-x_{3,2}y_{2,3}-y_{3,2}x_{2,3}.$
  \item $Q^{(3)}(A)=\det A$ and $\Theta_{Q^{(3)}}(X,Y,Z)$ is obtained from  $Q^{(3)}(A)$ replacing $a_{i,j}a_{k,l}a_{r,s}$ by the sum
  $x_{i,j}y_{k,l}z_{r,s}+x_{k,l}y_{i,j}z_{r,s}+x_{i,j}y_{r,s}z_{k,l}+x_{r,s}y_{k,l}z_{i,j}+x_{k,l}y_{r,s}z_{i,j}+x_{r,s}y_{i,j}z_{k,l}$
\end{enumerate}
\end{enumerate}
To find the expression of the bilinear form $\Theta_{Q^{(2)}}$ from the quadratic form $Q^{(2)}(A)$ is classical: we use the formula
$$4\Theta_{Q^{(2)}}(X,Y)=Q^{(2)}(X+Y)-Q^{(2)}(X-Y).$$
We have similar relations for greater degree, for example
$$24\Theta_{Q^{(3)}}(X,Y,Z)=Q^{(3)}(X+Y+Z)-Q^{(3)}(X+Y-Z)-Q^{(3)}(X-Y+Z)-Q^{(3)}(-X+Y+Z).$$

\medskip

Let $A$ be a matrix of $gl(n,\K)$ and $\widetilde{A}$ be a $^*\K$-perturbation of $A$ that is $\widetilde{A}=A+E$ where $E$ is a matrix of $gl(n,I)$ that is $E=(\epsilon_{i,j})$ with $\epsilon_{i,j} \in I$. The characteristic polynomial of $\widetilde{A}$ is $^*\K$-perturbation of $C_A(X)$ implying
$$C_{\widetilde{A}}(X)=C_A(X)+\Xi(X)$$
with $\Xi(X) \in I[X].$
We have
$$C_{A+E}(X)= X^n-Q^{(1)}(A+E)X^{n-1}+ \cdots + (-1)^jQ^{(j)}(A+E)X^{n-j}+\cdots+(-1)^nQ^{(n)}(A+E)$$
with
\begin{enumerate}
  \item $Q^{(1)}(A+E)={\rm tr}(A+E)=Q^{(1)}(A)+Q^{(1)}(E)$
  \item $Q^{(2)}(A+E)=Q^{(2)}(A)+ Q^{(2)}(E)+\Theta_{Q^{(2)}}(A,E)$
  \item $Q^{(3)}(A+E)=Q^{(3)}(A)+ Q^{(3)}(E)+\ds\frac{1}{2}\Theta_{Q^{(3)}}(A,A,E)+\ds\frac{1}{2}\Theta_{Q^{(3)}}(A,E,E)$
  \item  $Q^{(4)}(A+E)=Q^{(4)}(A)+ Q^{(4)}(E)+\ds\frac{1}{6}\Theta_{Q^{(4)}}(A,A,A,E)+\ds\frac{1}{4}\Theta_{Q^{(4)}}(A,A,E,E)$
   
   $+\ds\frac{1}{6}\Theta_{Q^{(4)}}(A,E,E,E)$
\end{enumerate}
and more generally
$$\begin{array}{ll}
Q^{(k)}(A+E)=      &  Q^{(k)}(A)+\ds\frac{1}{(k-1)!}\Theta_{Q^{(k)}}(A,\cdots,A,E)+ \ds\frac{1}{2!(k-2)!}\Theta_{Q^{(k)}}(A\cdots,A,E,E)  \\
&  + \cdots +\ds\frac{1}{i!(k-i)!}\Theta_{Q^{(k)}}(A,\cdots,A,E,\cdots,E)
\\
&+\cdots +\ds\frac{1}{(k-1)!}\Theta_{Q^{(k)}}(A,E,\ \cdots,E)+ Q^{(k)}(E).
\end{array}
$$
We deduce
\begin{proposition}
$$\begin{array}{ll}
\Xi(X)=&-Q^{(1)}(E)X^{n-1}+(Q^{(2)}(A+E)-Q^{(2)}(A))X^{n-2}+ \cdots \\
&+(-1)^{n-k}Q^{(k)}(A+E)-Q^{(k)}(A)X^k+ \cdots+(-1)^n Q^{(n)}(A+E)-Q^{(n)}(A))).
\end{array}$$
\end{proposition}
We can give the approximation to the first order (the notion is obvious) of $C_{\widetilde{A}}(X)$. If $E=(\varepsilon_{i,j})$ with $\varepsilon_{i,j} \in I$, its decomposition is
$$E=\alpha_1U_1+\alpha_1\alpha_2U_2+ \cdots +\alpha_1\alpha_2\cdots \alpha_{n^2}U_{n^2}$$ with $U_1,U_2,\cdots,U_{n^2}$
linearly independent in $gl(n,\K)$. The approximation of the first order is given by $\alpha_1$. Since
$$Q^{(k)}(A+E)-Q^{(k)}(A)\simeq \ds\frac{1}{(k-1)!}\Theta_{Q^{(k)}}(A,\cdots,A,E)$$
the approximation corresponds to
$$Q^{(k)}(A+E)-Q^{(k)}(A)\simeq \ds\frac{1}{(k-1)!}\Theta_{Q^{(k)}}(A,\cdots,A,\alpha_1U_1).$$
\begin{theorem}\label{xi2}
Let $A$ be a matrix of $gl(n,\K)$ and $\widetilde{A}$ be a $^*\K$-perturbation of $A$ that is $\widetilde{A}=A+E$ where $E=(\varepsilon_{i,j})$ is a matrix of $gl(n,I)$. The expression of the characteristic polynomial $C_{\widetilde{A}}(X)$ of $\widetilde{A}$ which is a $^*\K$-perturbation of $C_A(X)$ is
$$C_{\widetilde{A}}(X)=X^n-Q^{(1)}(A+E)X^{n-1}+ \cdots + (-1)^jQ^{(j)}(A+E)X^{n-j}+\cdots+(-1)^nQ^{(n)}(A+E)$$
where
$$Q^{(k)}(A+E)=  Q^{(k)}(A)+\sum_{i=1}^{k} \ds\frac{1}{i!(k-i)!}\Theta_{Q^{(k)}}(\underbrace{A,\cdots,A}_{k-i},\underbrace{E,\cdots,E}_{i})$$
and $\Theta_{Q^{(k)}}(A_1,\cdots,A_k)$ the symmetric $k$-linear form associated with the degree $k$ homogeneous polynomial  $Q^{(k)}(A)$.
If 
$$E=(\varepsilon_{i,j})=\alpha_1U_1+\alpha_1\alpha_2 U_2 +\cdots +\alpha_1\alpha_2\cdots\alpha_lU_l$$ is the decomposition of $E$ with 
$U_s=(u_{i,j}^{(s )})\in gl(n,\K)$ then the approximation at the order $1$ of $C_{\widetilde{A}}(X)-C_A(X)$ is given by the polynomial
$$\Xi_1(X)= \alpha_1\sum_{k=1}^{n}(-1)^k\ds\frac{1}{(k-1)!}\Theta_{Q^{(k)}}(A,\cdots,A,U_1)X^{n-k}.$$
\end{theorem}

\subsection{Linear deformations}
We call linear deformation of a matrix $A \in gl(n,\K)$ a  $^*\K$-perturbation of $A $ of the form
$$\widetilde{A}=A+ \varepsilon U_0$$ with  $U_0 \in gl(n,\K)$ and $\varepsilon\in \mathcal{I}$.  From Theorem \ref{xi2},
$$C_{\widetilde{A}}(X)=C_A(X)+\Xi(X)$$
with
$$
\begin{array}{l}
     \Xi(X)=  \varepsilon\left(\Theta_{Q^{(1)}}(U_0)X^{n-1}-\Theta_{Q^{(2)}}(A,U_0)X^{n-2}+\cdots +(-1)^{n-1}\Theta_{Q^{(n)}}(A,A,\cdots,A,U_0) \right)  \\
       +\varepsilon^2(-\Theta_{Q^{(2)}}(U_0,U_0)X^{n-2}+\Theta_{Q^{(3)}}(A,U_0,U_0)X^{n-3}+ \cdots +(-1)^{n-1}\Theta_{Q^{(n)}}(A,\cdots,A,U_0,U_0) ) \\
      
 \cdots  +\varepsilon^n (-1)^{n-1}\Theta_{Q^{(n)}}(U_0,U_0,\cdots,U_0).
\\
\end{array}
$$
\begin{proposition}
Let $\widetilde{A}=A+ \varepsilon U_0$ be a linear perturbation of $A \in gl(n,\K)$. Then the characteristic polynomial $C_{\widetilde{A}}(X)$ is a perturbation of $C_A(X)$ and
$$\frac{C_{\widetilde{A}}(X)-C_A(X)}{\varepsilon}\simeq \Theta_{Q^{(1)}}(U_0)X^{n-1}-\Theta_{Q^{(2)}}(A,U_0)X^{n-2}+\cdots (-1)^{n-1}\Theta_{Q^{(n)}}(A,A,\cdots,A,U_0) . $$
\end{proposition}
In particular $\ds \frac{C_{\widetilde{A}}(X)-C_A(X)}{\varepsilon} \in \mathcal{I}(X)$ if and only if
$$U_0 \in sl(n,\K) \ {\rm and} \  \Theta_{Q^{(k)}}(A,A, \cdots,A,U_0)=0, \ k=2, \cdots, n.$$
\begin{corollary}
We have $C_{\widetilde{A}}(X)=C_A(X)$ if and only if 
$$\Theta_{Q^{(k)}}(A,\cdots,A,U_0)=0$$
for any $k \geq 1$.
\end{corollary}
Assume now that $a$ is a simple root of $C_A(X)$ that is an eigenvalue of $A$. Then there exists $\xi \in \mathcal{I}$ such that $a+\xi$ is an eigenvalue of $\widetilde{A}$.  We have
$$C_{\widetilde{A}}(a+\xi)\simeq \xi C'_A(a)+ \varepsilon \Pi_1(a)$$
where $\Pi_1(X)=\Theta_{Q^{(1)}}(U_0)X^{n-1}-\Theta_{Q^{(2)}}(A,U_0)X^{n-2}+\cdots (-1)^{n-1}\Theta_{Q^{(n)}}(A,A,\cdots,A,U_0) $
as soon as $\Pi_1(a) \neq 0.$ In this case
$$\xi \simeq -\varepsilon \frac{\Pi_1(a)}{C'_A(a)}.$$
If $a$ is a root of order $p$, then
$$\xi ^p \simeq -\varepsilon \frac{p!\Pi_1(a)}{C^{(p)}_A(a)}.$$
\noindent{\bf Example.} Let $A$ be the matrix
$$\begin{pmatrix}
     0 & 1  & 0   \\
     0 & 0 & 1 \\
     0 & 0 & 0 
\end{pmatrix}
$$ 
and
$$\widetilde{A}=\begin{pmatrix}
     0 & 0  & 0   \\
     0 & 0 & 0 \\
     \varepsilon & 0 & 0 
\end{pmatrix}
$$ 
In this case $a=0$ is an eigenvalue of multiplicity $3$, $U_0$ is the matrix
$$U_0=\begin{pmatrix}
     0 & 0  & 0   \\
     0 & 0 & 0 \\
    1 & 0 & 0 
\end{pmatrix}
$$ 
and we have
$$\Theta_{Q^{(1)}}(U_0)=0, \Theta_{Q^{(2)}}(A,U_0)=0, \ \Theta_{Q^{(3)}}(A,A,U_0)=1$$
and $\Pi_1(X)=1.$ We deduce
$$\xi ^3 \simeq -\varepsilon.$$
\subsection{Conservative perturbations}
Let $A$ be in $gl(n,\K)$ and $\widetilde{A}=A+E$ a $^*\K$-perturbation of $A$. This perturbation is called conservative if the eigenvalues of $\widetilde{A}$ are those of $A$. This is equivalent to say that 
$$C_A(X)=C_{\widetilde{A}}(X)$$ 
or equivalently$$\Xi(X)=0.$$
This is equivalent to the system
$$
\left\{
\begin{array}{l}
tr (E)=0,\\
\Theta_{Q^{(2)}}(A,E)+Q^{(2)}(E)=0,\\

\cdots\\

\Theta_{Q^{(n)}}(\underbrace{A,\cdots,A}_{n-1},E)+\cdots+
\Theta_{Q^{(n)}}(\underbrace{A,\cdots,A}_{n-j},\underbrace{E,\cdots,E}_{j})+\cdots\\
\qquad +\Theta_{Q^{(n)}}(A,\underbrace{E,\cdots,E}_{n-1})+Q^{(n)}(E)=0.
\end{array}
\right.
$$
Thus, for a given matrix $A \in gl(n,\K)$, the conservative perturbations correspond to points of an algebraic variety $V_A$ contained in $sl(n,^*\K))$. 

\noindent{\bf Example: $n=2$.} The conservative perturbations of a given matrix $\ds A=\begin{pmatrix}
    a_1  &  a_2  \\
  a_3    &  a_4
\end{pmatrix}$
are given by the matrix $\ds E=\begin{pmatrix}
   \varepsilon_1  &   \varepsilon_2  \\
  \varepsilon_3    &  - \varepsilon_1
\end{pmatrix}$
satisfying
$$-(a_1-a_4)\varepsilon_1-a_2\varepsilon_3-a_3\varepsilon_2-(\varepsilon_1)^2-\varepsilon_2\varepsilon_3=0.$$
Using the decomposition $(\varepsilon_1,\varepsilon_2,\varepsilon_3)=\alpha_1U_1+\alpha_1\alpha_2U_2+\alpha_1\alpha_2\alpha_3U_3$
we obtain
$\alpha_1(-(a_1-a_4)u_{11}-a_2u_{31}-a_3u_{21}+\alpha_1\rho =0$ with $\rho \in I$ where $U_i=(u_{ji})$. This gives
$$-(a_1-a_4)u_{11}-a_2u_{31}-a_3u_{21}=0.$$
In particular, if the decomposition is of length $1$ that is $\alpha_2=0$, then
$$
\left\{
\begin{array}{l}
     -(a_1-a_4)u_{11}-a_2u_{31}-a_3u_{21} =0 \\
     -u_{11}^2-u_{21}u_{31}=0
\end{array}
\right.
$$
that is the matrix $$\widehat{U_1}
 =\begin{pmatrix}
    u_{11}  &  u_{21} \\
     u_{31} &   -u_{11}
\end{pmatrix}
$$
is singular and orthogonal to $A$ with respect the inner product $Q^{(2)}$ in $gl(n,\K)$.  For example, if $A=0$ then $V_A$ is the algebraic variety contained in $sl(2,^*\K)$ whose elements are the singular matrices.

\subsection{Particular case: Hermitian matrices}
Recall that  a matrix $A \in gl(n,\C)$ is called hermitian if it satisfies
$$^t \overline{A }=A.$$
The coefficients of the characteristic polynomial of an hermitian matrix are real and the eigenvalues are also real. Then we can consider that these eigenvalues are ordered:
$$\lambda_1 \leq \lambda_2 \leq \cdots \leq \lambda_n.$$
Let $\widetilde{A}=A+E$ be an hermitian $^*\C$-perturbation of $A$, that is $\widetilde{A}$ is an hermitian matrix. This implies that also $E$ is an hermitian matrix in $gl(n,I)$. The eigenvalues $\widetilde{\lambda_i} $   of $\widetilde{A}$ are in  $^*\R$ and we have 
$$\widetilde{\lambda_1} \leq \widetilde{\lambda_2 }\leq \cdots \leq \widetilde{\lambda_n}.$$
We can write
$$E=\alpha_1U_1+\alpha_1\alpha_2U_2+\cdots+\alpha_1\alpha_2\cdots\alpha_pU_p$$
where the matrices $U_i$ are hermitian and linearly independent in $gl(n,\C)$,  
$\alpha_i \in I, p\leq n^2.$ Let $C_{\widetilde{A}}(X) $ and $C_A(X)$  be the characteristic polynomials of $\widetilde{A}$ and $A$. We have $C_{\widetilde{A}}(X) =C_A(X)+\Xi (X)$ with $\Xi(X) \in \mathcal{I}I[X]$. Let $\widetilde{\lambda}=\lambda+\rho $ be an eigenvalue of $\widetilde{A}$, with $\lambda$ an   eigenvalue of $A$ and $\rho \in \mathcal{I}$. If $\lambda$ is a simple root, 
$$C_{\widetilde{A}}(\lambda+\rho)=0= C_A(\lambda+\rho)+\Xi(\lambda+\rho)\simeq \rho C_A'(\lambda)+\Xi_1(\lambda)$$
with 
$$\Xi_1(\lambda)= \alpha_1\sum_{k=1}^{n}(-1)^k\ds\frac{1}{(k-1)!}\Theta_{Q^{(k)}}(A,\cdots,A,E)\lambda^{n-k}.$$
For example, for $n=2$ we have
\begin{proposition}
If $\lambda$ is a simple eigenvalue of $A\in \mathcal{M}_2(\mathbb{C})$ then
$$\rho \simeq \alpha_1 \ds \frac{-\lambda tr( U_1)+\theta_{Q^{2}}(A,U_1)}{C'_A(\lambda)}$$
as soon as $-\lambda tr( U_1)+\Theta_{Q^{2}}(A,U_1) \neq 0.$
\end{proposition}
\medskip

If $\lambda$ is the largest eigenvalue of $A$, since $A$ is hermitian, we have
$$\lambda=\max \ ^t\overline{X}AX$$
for $X \in \K^n$ and $^t\overline{X}X=1$.  If $\widetilde{A}=A+E$ is an hermitian $^*\C$-perturbation of $A$, then the largest eigenvalue $\widetilde{\lambda}$ of $\widetilde{A}$ satisfies
$$\widetilde{\lambda} =\max \ ^t\overline{\widetilde{X}}\widetilde{A}\widetilde{X}$$
for $\widetilde{X} \in \, ^ *\C^n$ and $^t\overline{\widetilde{X}}X=1$. But we have seen that $\widetilde{\lambda} =\lambda + \rho$ with $\rho \in I$.  This last condition implies that the vector $\widetilde{X}$ which belongs to $^*\C^n$ has no infinitely large component, that is whose inverse is in $I$. This implies that $\widetilde{X}=X+Y$ with $X \in \C^n$ and $Y \in \mathcal{I}^n$. We deduce
$$1=\, ^t\overline{\widetilde{X}}X=\, ^t\overline{X}X+\, ^t\overline{X}Y+X^t\overline{Y}+\, ^t\overline{Y}Y.$$
But all the terms are infinitesimal except $^t\overline{X}X$ implying
$$^t\overline{X}X=1, \ ^t\overline{X}Y+X^t\overline{Y}+\, ^t\overline{Y}Y=0$$
and also
$$^t\overline{X}Y+X^t\overline{Y}=0, ^t\overline{Y}Y=0.$$
Since the inner product $^t\overline{U}U$ is non degenerate, we obtain $Y=0$ and $\widetilde{X}=X$. Then the largest eignevalue of $\widetilde{A}$ is given by the formula
$$\widetilde{\lambda} =\max \ ^t\overline{X}\widetilde{A} X$$
for $X \in \C^n$ and $^t\overline{X}X=1$. So
$$\widetilde{\lambda} =\max \ (^t\overline{X}A X+\, ^t\overline{X}E X).$$

\end{document}